\documentclass[10pt]{article}

\usepackage{amsmath}
\usepackage{amsfonts}

\newcommand{\be}{\begin}
\newcommand{\Ref}[1]{(\ref{#1})}

\newcommand{\ca}{\mathcal A}
\newcommand{\cb}{\mathcal B}
\newcommand{\cg}{\mathcal G}

\newcommand{\al}{\alpha}
\newcommand{\del}{\delta}
\newcommand{\eps}{\epsilon}
\def\ga{\gamma}

\def\ka{\kappa}
\def\lla{\lambda}
\def\La{\Lambda}

\def\Om{\Omega}
\newcommand{\si}{\sigma}
\def\Si{\Sigma}

\newcommand{\R}{\mathbb R}
\newcommand{\z}{\mathbb Z}
\def\q{{\mathbb Q}}
\newcommand{\co}{\mathbb C}
\newcommand{\e}{{\mathbb E}}
\def\f{{\mathbb F}}

\def\ba{{\mathbb A}}

\def\ub{{\mathbb U}}

\newcommand{\oset}{\overset}
\newcommand{\uline}{\underline}
\newcommand{\oline}{\overline}

\newcommand{\la}{\langle}
\newcommand{\ra}{\rangle}
\newcommand{\st}{\;|\;}
\newcommand{\ti}{\tilde}
\newcommand{\setmin}{\backslash}
\newcommand{\prtl}{\partial}
\def\square{\kern20pt{\vbox{\hrule height.4pt
	\hbox{\vrule width.4pt height 6pt\kern6pt
		\vrule width.4pt}
	\hrule height.4pt}}}
\newcommand{\ry}{\R\times Y}
\newcommand{\rpy}{\R_+\times Y}

\def\im{\text{im}}

\newcommand{\fact}{\frac1{8\pi^2}}
\def\inv{^{-1}}

\def\proof{{\em Proof.}\ }

\def\frakg{\mathfrak g}
\def\cs{\text{cs}}
\def\csp{\text{cs}_\pi}

\def\U#1{\text{U}(#1)}
\def\SU#1{\text{SU}(#1)}
\def\su{\SU2}
\def\SO#1{\text{SO}(#1)}
\def\so{\SO3}
\def\hf{\text{HF}}
\def\hfh{\widehat{\hf}\vphantom{\hf}}

\def\poinc{\text{Poincar\'e}}

\def\cem{\check{M}}
\def\harm{\mathcal H}
\def\cl{\mathcal L}
\def\ct{\mathcal T}

\def\torsion{\text{torsion}}

\def\bz{{\mathbf z}}

\def\cbb{\cb^*((-1,1)\times Y)}
\def\band{(-1,1)\times Y}
\def\cf{\text{CF}}

\begin{document}

\title{Equivariant aspects of Yang-Mills Floer theory}

\author{Kim A.\ Fr\o yshov\thanks{Partially supported by
a post-doctoral grant from the Norwegian Research Council,
and by NSF grant DMS-9971731.}}

\maketitle

\bibliographystyle{plain}

\section{Introduction}

This paper is concerned with Floer cohomology groups
of $\so$ bundles $P\to Y$, where $Y$ is a closed, oriented
$3$-manifold. Following \cite{BD1} we only consider
{\em admissible} bundles, which means $P$ should be non-trivial
over some surface in $Y$ unless $Y$ is an (integral) homology sphere.
The mod~$8$ periodic Floer group $\hf^*(P;G)$ with coefficients in the abelian
group $G$ is then a topological invariant of $Y,P$. We will often omit
the coefficient group from notation. When nothing else is
specified our results hold for any coefficient group.

If $P$ is non-trivial then ``cup product'' with a certain
$4$-dimensional cohomology class (four times the $\mu$-class of a point)
defines a homomorphism $u:\hf^*(P)\to
\hf^{*+4}(P)$. We use Floer's exact triangle to show that there is
always a positive integer $n$ such that
\be{equation*}
(u^2-64)^n=0.
\end{equation*}

When $Y$ is a homology sphere the $u$-map is in general defined on $\hf^q(Y)$
only for $q\not\equiv4,5\mod8$, due to the presence of the trivial $\su$
connection over $Y$. However, by ``factoring out'' interaction with the
trivial connection we construct a {\em reduced} Floer group $\hfh^*(Y)$
on which the $u$-map is defined in all degrees. We expect that
$\hfh^*(Y;\R)$ is isomorphic to the $\R[u]$ torsion submodule
of Austin and Braam's
equivariant Floer group \cite{AuB2}, but this will not be proved in
this paper.

For the reduced Floer group we again find that $(u^2-64)^n=0$ for some 
$n$. Combining this with the splitting theorem of
Freedman and Taylor \cite{FT}, \cite{Stong} we obtain a proof,
in the simply-connected case, of the finite
type conjecture of Kronheimer and Mrowka:
\be{thm}Let $X$ be a smooth, compact, simply-connected, oriented
$4$-manifold with $b_2^+(X)$ odd and $\ge3$. Then the Donaldson invariants
of $X$ have finite type.
\end{thm}
See Section~\ref{fintypsec} for more details.
Quite different proofs of this theorem have been given by Mu\~noz~\cite{mun2}
(without any assumptions on the fundamental group) and by
Wieczorek~\cite{wiec1}.

The nilpotency of $u^2-64$ also leads to the following theorem, by another
application of Floer's exact triangle:
\be{thm}If $Y$ is any oriented homology $3$-sphere and $R$ any
associative ring in which
$2$ is invertible then $u:\hf^q(Y;R)\to\hf^{q+4}(Y;R)$ is an isomorphism
for $q\not\equiv4,5\mod8$. In particular, $\hf^*(Y;R)$ is mod~$4$ periodic.
\end{thm}

We now focus on rational coefficients. Let $Y$ be an oriented homology
$3$-sphere.
Computing the reduced Floer group $\hfh^*(Y;\q)$ from the ordinary
one $\hf^*(Y;\q)$ merely requires the knowledge
of a single integer $h(Y)$,
which measures interaction between irreducible
flat $\su$ connections and the trivial connection over $Y$.
This invariant $h(Y)$ has the following properties:
\be{thm}\label{h-prop}\
\be{description}
\item[(i)]$h(Y_1\#Y_2)=h(Y_1)+h(Y_2)$.
\item[(ii)]If the homology sphere $Y$ bounds a
smooth, compact oriented $4$-manifold
with negative definite intersection form then $h(Y)\ge0$, with strict 
inequality if the intersection form is not diagonal over the integers.
\item[(iii)]For the Brieskorn sphere $\Si(2,3,5)$ one has $h=1$.
\end{description}
\end{thm}
Since $h(S^3)=0$, property (ii) generalizes Donaldson's theorem \cite{D1},
\cite{D2}.

Theorem~\ref{h-prop} implies that
$h$ is a surjective group homomorphism
\[h:\theta^H_3\to\z,\]
where $\theta^H_3$ is the integral homology
cobordism group of oriented homology $3$-spheres.
The mod~$2$ reduction of $h$ is not the
Rochlin invariant, because $\Si(2,3,7)$ bounds
an orientable rational ball but has Rochlin invariant one \cite{FS6}.

Let $k$ be a positive integer and $\ga$ a negative knot in $S^3$
(ie a knot which admits a regular projection with only negative crossings).
The manifold $S^3_{\ga,1/k}$ resulting from $1/k$ surgery on $\ga$ bounds
both positive and negative definite $4$-manifolds, hence $h(S^3_{\ga,1/k})=0$.
On the other hand, if $\ga$ is non-trivial then, as proved in \cite{CG},
$S^3_{\ga,-1/k}$ bounds
a $4$-manifold with negative definite intersection form of the type
$-E_8\oplus n(-1)$ for 
some $n\ge0$, hence $h(S^3_{\ga,-1/k})>0$.
For instance, if $p,q$ are mutually prime integers $\ge2$
then the Brieskorn sphere $\Si(p,q,pqk\pm1)$ results
from $\pm1/k$ surgery on the negative $(p,q)$ torus knot, so
\[h(\Si(p,q,pqk+1))=0,\qquad h(\Si(p,q,pqk-1))>0.\]

Further properties of the $h$-invariant will be given in
\cite{Fr6}, including estimates on the behaviour under surgery on knots.

If one takes a field of characteristic $p>2$ as coefficients for the
Floer groups then one obtains a group homomorphism
\[h_p:\theta^H_3\to\z\]
(which depends only on $p$) with the following property: if $Y$ bounds 
a negative definite, smooth, compact, oriented $4$-manifold
without $p$-torsion in its homology then $h_p(Y)\ge0$,
with strict inequality if the intersection form is non-standard.
Unfortunately, the author is unable to prove anything
about $h_p$ that does not also hold for $h$. The invariants $h_p$ will
therefore not be pursued any further in this paper.

We intend to discuss related topics in Seiberg-Witten Floer
theory, and in Yang-Mills Floer theory with $\z/2$ coefficients, in
forthcoming papers.

\section{Preliminaries}\label{Floer-homology}

This section is mostly a review of well known material.
For more details see \cite{DK}, \cite{F1},
\cite{BD1}, \cite{DS1}, \cite{D5}.

If $X$ is an $n$-dimensional smooth manifold, with or without boundary,
and $E\to X$ a rank~$2$ 
unitary vector bundle let $\ca(E)$ denote the space of all connections in 
$E$ which induce a fixed connection in $\La^2E$. Let $\cg(E)$ be the
group of automorphisms (or gauge transformations)
of $E$ of determinant~$1$ and set
$\cb(E)=\ca(E)/\cg(E)$.

In the case of an $\so$ bundle $P\to X$ we define $\cg(P)$ to be the group of
all automorphisms of $P$ and $\cg_S(P)$ to be the subgroup of {\em even} 
automorphisms, ie automorphisms that lift to sections of
$P\times_{Ad(\so)}\su$. Note that there is a natural exact sequence
\be{equation}\label{gauge:ex}
1\to\cg_S(P)\to\cg(P)\oset{\eta}\to H^1(X;\z/2)\to0.
\end{equation}
If $\frakg_E$ is the $\so$ bundle associated to the $U(2)$ bundle $E$ then
we can identify $\cb(E)$ with the space of all connections in $\frakg_E$ 
modulo even automorphisms.
A connection in $E$ is called {\em irreducible} if its stabilizer in 
$\cg(E)$ is $\{\pm1\}$; otherwise it is called {\em reducible}. We say 
a connection in $E$ is {\em twisted reducible} if the induced connection in 
$\frakg_E$ respects a splitting $\frakg_E=\lla\oplus L$, where $\lla$ is 
a non-trivial real line bundle and $L$ a non-orientable real $2$-plane
bundle.

If $X$ is compact and we consider $L^p_1$ connections in $E$ modulo
$L^p_2$ gauge transformations, where $p>n/2$, then the subspace
$\cb^*(E)\subset\cb(E)$ of irreducible connections is a Banach 
manifold; if in addition $p$ is an even integer then $\cb^*(E)$ admits 
$C^\infty$ partitions of unity.

\subsection{Floer cohomology groups of homology $3$-spheres}

Let $Y$ be an oriented (integral) homology $3$-sphere. The Floer
cohomology group $\hf^*(Y)$ was defined in \cite{F1} by applying Morse
theoretic ideas to a suitably perturbed Chern-Simons function
$\csp:\cb(Y\times\su)\to\R/\z$. A critical point $A$ of $\cs_\pi$ is called
{\em non-degenerate} if the Hessian is non-singular
on $\ker(d^*_A)\subset\Om^1_Y(su(2))$. For the unperturbed Chern-Simons 
function $\cs$ the critical points are the flat connections and the
Hessian is $*d_A$. In this case  
there is up to gauge equivalence only one reducible critical point, namely
 the trivial connection $\theta$, which is non-degenerate since 
$H^1(Y;\R)=0$. The trivial connection is in fact a critical point of
$\csp$ for all perturbations $\pi$ (this is a consequence of the gauge
invariance of $\csp$). We will always work with a small, generic perturbation;
then we may assume that $\csp$ has only
finitely many critical points, all of which
are non-degenerate, and that $\theta$ is the only reducible 
critical point. 
We will usually not refer explicitly to the perturbation $\pi$ or the
corresponding perturbations of the anti-self-dual equations.

For any pair of flat $\su$ connections $\al,\beta$ over $Y$ and any
real number $\ka$ with $\ka\equiv\cs(\al)-\cs(\beta)\mod\z$ let
$M(\al,\beta;\ka)$
be the moduli space of all anti-self-dual $\su$ connections $A$ over $\ry$ which are
asymptotic to $\al$ at $\infty$ and to $\beta$ at $-\infty$,
and have second relative Chern class $\fact\int_{\ry}\text{tr}(F_A^2)=\ka$.
Here $F_A$ is the curvature of $A$.
These moduli spaces are orientable, and orientations should be chosen 
compatible with gluing maps and addition of instantons over
$S^4$, see \cite{D2}. We denote by $M(\al,\beta)$ the
moduli space $M(\al,\beta;\ka)$ whose expected dimension lies in the interval
$[0,7]$, and set
\[\cem(\al,\beta)=M(\al,\beta)/\R,\]
where $\R$ acts by
translation. For an irreducible flat connection $\al$ we define the index
$i(\al)\in\z/8$ by
\[i(\al)\equiv\dim M(\al,\theta)\mod8.\]

The Floer cohomology group $\hf^*(Y)$ is the cohomology of the $\z/8$ graded
cochain complex $(\cf^*,d)$, where $\cf^i$ is the free abelien group
generated by the gauge equivalence classes of
irreducible flat $\su$ connections of index $i$ over $Y$. The differential
$d$ has matrix coefficient $\#\,\cem(\al,\beta)$ when $i(\al)-i(\beta)=1$, 
where $\#$ means the number of points counted with sign. To show that
$d^2=0$ one counts the ends of $\cem(\al,\ga)$ when $i(\al)-i(\ga)=2$.

The Floer homology group $\hf_*(Y)$ is the homology of the dual complex
of $(\cf^*,d)$. 
There is then a canonical identification $\hf_q(Y)=\hf^{5-q}(\oline Y)$.

The trivial connection over $Y$ gives rise to a homomorphism
$\del:\cf^4\to\z$ and an element $\del'\in \cf^1$, defined by
\[\del\al=\#\cem(\theta,\al),\qquad
\del'=\sum_\beta\#\cem(\beta,\theta)\,\beta,\]
where $\beta$ runs through the generators of $\cf^1$.
These satisfy $\del d=0$ and $d\del'=0$ (for the same reason
that $d^2=0$) and so define
\be{equation*}
\del_0:\hf^4(Y)\to\z,\quad\del'_0\in\hf^1(Y).
\end{equation*}
These will play a central role in this paper.

\subsection{Floer cohomology groups of non-trivial $\so$ bundles}

Now let $Y$ be a closed, oriented $3$-manifold and $P\to Y$ a non-trivial
$\so$ bundle which is 
admissible in the sense of \cite{BD1}. This means that the Stiefel-Whitney
class $w_2(P)$ is not the
mod~$2$ reduction of a torsion class in $H^2(Y;\z)$, or equivalently that
$w_2(P)$ defines a non-zero map $H_2(Y;\z)\to\z/2$. In particular, the
Betti number $b_1(Y)$ must be positive. In this case there are no reducible
flat connections in $P$. The Floer group $\hf^*(P)$ is defined just as
for homology spheres, using a small, generic perturbation of the
Chern-Simons function. However, the group is now only {\em affinely}
$\z/8$ graded, ie only the index difference of two flat
connections is well defined in $\z/8$.

On the other hand, given a {\em spin structure} on $Y$ 
one can define a  mod~$4$ grading on $\hf^*(P)$ as follows.
Let $\al$ be a (non-degenerate) flat connection in $P$. As shown in 
\cite{Kirby2}, $Y$ spin bounds a simply-connected, spin $4$-manifold $X'$. 
Choose a cylindrical end metric on 
the corresponding open $4$-manifold $X$.
Since the adjoint vector bundle of $P$
is isomorphic to $\R\oplus L$ for some complex line bundle
$L\to Y$, and since the restriction map $H^2(X)\to H^2(Y)$ is
surjective, there is an $\so$ bundle $Q\to X$ with
$Q|_Y\approx P$. Now define the index $i(\al)\in\z/4$ by
\[i(\al)\equiv\dim M(Q,\al)+3b^+_2(X)\mod4.\]
Here $M(Q,\al)$ is the moduli space of anti-self-dual connections in the 
bundle $Q$ which are asymptotic to $\al$ over the end,
while $b^+_2(X)$ is the maximal dimension of a positive subspace
for the intersection form on $H_2(X;\q)$. It follows easily from
the dimension formula for anti-self-dual moduli spaces over closed
$4$-manifolds that $i(\al)$ is well defined mod~$4$.

Recall that the group $H^1(Y;\z/2)$ acts simply transitively on the set of
spin structures on $Y$. If $s_1,s_2$ is a pair of spin structures then the
corresponding index functions $i_1,i_2$ are related by
\[(i_1-i_2)/2\equiv
((s_1-s_2)\cup w_2(P))[Y]\mod2.\]
In particular, the grading mod~$2$ is independent of the spin structure.

To any class $\zeta\in H^1(Y;\z/2)$ we can associate an involution
$\zeta^*$ of $\hf^*(P)$, well defined up to an
overall sign; this involution is induced by the twisted bundle
\[(\R_-\times P)\cup_g(\R_+\times P)\]
over $\ry$ where $g$ is any automorphism of $P$ with $\eta(g)=\zeta$
(recall the exact sequence
\Ref{gauge:ex}). The reason for the sign ambiguity is that one has to make a
choice concerning the orientation of moduli spaces in this bundle.
In any case, one does get a group homomorphism
\[H^1(Y;\z/2)\to\text{Aut}(\hf^*(P))/\{\pm1\}.\]
It is easy to see that each $\zeta^*$ has degree $0$ or $4$, and if 
$\zeta$ has an integral lift then
\[\deg(\zeta^*)/4\equiv(w_2(P)\cup\zeta)[Y]\mod2.\]
(See \cite{BD1} for the general formula.) In particular, there is always 
a degree~$4$ involution $\zeta^*$.

\subsection{Invariants of $4$-manifolds with boundary}

The purpose of this subsection is merely to review the definition
of Donaldson invariants of $4$-manifolds with boundary. 
A gluing theorem for these invariants will be stated in
Section~\ref{fintypsec}.

Let $X$ be a smooth, oriented Riemannian $4$-manifold with one tubular end
$\rpy$, and let $E\to X$ be a $U(2)$ bundle. For simplicity, suppose
the $\so$ bundle $P\to Y$ associated to $E|_Y$ is non-trivial and admissible.
In this case the Donaldson invariant for the bundle $E$ is a linear map
\be{equation}\label{psie}
D_E:\ba(X)\to \hf^*(P;\q),
\end{equation}
well defined up to an overall sign, where
\[\ba(X)=\text{Sym}(H_{\text{even}}(X;\q))\otimes\La(H_{\text{odd}}(X;\q)).\]
(When $Y$ is a homology sphere care must be taken to handle reducible 
connections in the moduli spaces; insofar as the corresponding invariants
are defined in this case they will be denoted
$D^c_X$, where $c=c_1(E)$.)

This invariant is defined in much the same way as the instanton invariants of
closed $4$-manifolds (see \cite{D3}, \cite{KM3}).
Let $\Si_1,\dots,\Si_m\subset X$ be a collection of smooth, compact, connected
submanifolds without boundary and in general position. As in \cite{KM3},
Section~2~(ii)
choose for $j=1,\dots,m$ a smooth, compact, codimension $0$ submanifolds
$U_j$ of $X$ containing $\Si_j$ such that the map
$H_1(U_j;\z/2)\to H_1(X;\z/2)$ is surjective, and such that the sets
$U_1,\dots,U_m$ are disjoint.
Let $\cb^*(U_j)=\cb^*(E|_{U_j})$ be
as in Section~\ref{Floer-homology} and let
$\e_j\to\cb^*(U_j)\times U_j$ be the universal $\so$ bundle. As in \cite{KM3}
choose a generic geometric representative
$V_j\subset\cb^*(U_j)$ for the cohomology class
$\mu(\Si_j)=-\frac14p_1(\e_j)/[\Si_j]$.
Let $d=\sum_j(4-\dim \Si_j)$ and for any flat connection $\al$ in $P$ set
\[Z_\al=r_1\inv(V_1)\cap\cdots\cap r_m\inv(V_m)\subset M(E,\al).\]
Here $M(E,\al)$ is the moduli space of (projectively) anti-self-dual connections in $E$
which are asymptotic to $\al$ at the end, and
$r_j:M(E,\al)\to\cb^*(U_j)$ is the restriction map.
For the monomial $z=[\Si_1]\cdots[\Si_m]\in\ba(X)$ define
\[D_E(z)=[\sum_\al(\#Z_\al)\al]\in \hf^*(P;\q),\]
where $\al$ runs through the equivalence classes of flat connections in $P$
for which $Z_\al$ has dimension $0$.

To show that $D_E(z)$ is independent of $U_j,V_j$ and
linear one can follow the arguments in
\cite{D3} and \cite{KM3} and show that when computing $D_E(z)$ one of 
the classes $\mu(\Si_j)$ may be evaluated ``abstractly''
(eg using \v Cech-type (co)homology).

\section{The $u$-map and the reduced Floer group}\label{umap}

In this section one could use any coefficient group for the Floer 
cohomology groups,
but for simplicity we will work with integral coefficients.

If $P$ is a non-trivial admissible $\so$ bundle over a closed, oriented
$3$-manifold $Y$, then the $u$-map $\hf^*(P)\to \hf^{*+4}(P)$ is defined,
roughly speaking, by evaluating the $4$-dimensional class $4\mu(x)$
over $4$-dimensional moduli spaces
$M(\al,\beta)$, where $\al,\beta$ are flat connections in $P$. If
$Y$ is a homology sphere, then the construction can still be carried out
on cochain level to give a homomorphism $v:\cf^*(Y)\to \cf^{*+4}(Y)$ (which
depends on certain choices). But due to the presence of the trivial connection,
this homomorphism is not quite a cochain map
(see Theorem~\ref{chain-map} below), and in general only defines a homomorphism
$u:\hf^q(Y)\to \hf^{q+4}(Y)$ for $q\not\equiv4,5\mod8$. However, by 
``factoring out'' interaction with the trivial connection we will 
construct a {\em reduced} Floer group $\hfh^*(Y)$, in which the $u$-map
is defined in all degrees.

\subsection{The $u$-map}

Let $Y$ be a closed, oriented $3$-manifold and $P\to Y$ an admissible
$\so$ bundle. We will define a graded homomorphism
$v:\cf^*(P)\to \cf^{*+4}(P)$. Let $\al$, $\beta$ be flat connections in $P$,
not both reducible.
Let $\e=\e(\beta,\al)\to M(\beta,\al)$ and
$\f\to\cb^*((-1,1)\times P)$ be the natural
oriented, euclidean $3$-plane bundles associated to the base-point
$(0,y_0)$. Here $(-1,1)\times P$ is the obvious $\so$ bundle.
There is a natural restriction map $r:M(\beta,\al)\to\cb^*((-1,1)\times P)$,
and we have $r^*\f=\e$. 
Choose sections $s_1,s_2$ of the complexified bundle
$\f\otimes\co$ and let $\si_j=r^*s_j$ be the induced sections of
$\e\otimes\co$. If $\dim\,M(\beta,\al)\le5$ then after perturbing
the $s_j$'s we may assume $\si_1$ has no zeros and that the section
$\si=\si_2\mod\si_1$ of the quotient bundle
$(\e\otimes\co)/\co\si_1$ is transverse to the
zero section. If $\al,\beta$ are both irreducible and
$\dim\,M(\beta,\al)=4$ then $\si\inv(0)$
is a finite set of oriented points, and we define the matrix coefficient
$\langle v(\al),\beta\rangle$ by
\[\langle v(\al),\beta\rangle=\#\si\inv(0).\]

The following theorem is due to Donaldson and Furuta \cite{D5}, but
we include a proof for the
sake of completeness and because we will need certain generalizations
later.

\be{thm}[Donaldson and Furuta]\label{chain-map}\
\be{description}
\item[(i)]If $P$ is a non-trivial, admissible bundle then $dv-vd=0$.
\item[(ii)]If $Y$ is a homology $3$-sphere then
\[dv-vd+2\del\otimes\del'=0,\]
where by
definition $\del=0$ in degrees $\not\equiv4\mod8$.
\end{description}
\end{thm}

In case (i) it follows that $v$ induces a homomorphism 
$u:\hf^*(P)\to \hf^{*+4}(P)$, while in case~(ii) one gets $u$-maps
\be{gather*}
\hf^i(Y)\to \hf^{i+4}(Y),\qquad i\neq4,5\\
\ker(\del_0)\to \hf^0(Y);\qquad \hf^5(Y)\to \hf^1(Y)/(\z\del'_0).
\end{gather*}

\proof Since (i) is essentially a special case of (ii),
we focus on the latter. Let $\al$, $\beta$ be irreducible, flat
$\su$ connections over a homology $3$-sphere $Y$ such that
$i(\beta)\equiv i(\al)+5\mod8$. We will show that
\[\langle(dv-vd+2\del\otimes\del')\al,\beta\rangle=0.\]
Our plan is to modify the section $\si$ for 
connections that are close to the trivial connection over $\band$, in order to
gain control over the ends of $\si\inv(0)$. Counting the number of
such ends with sign (this number must be zero) will then give (ii).

The proof is divided into four parts.

{\bf(I)} We first show that suitable modifications can be made to
the sections $s_j$ without affecting the definition of the chain map $u$.
This will be used in (IV) (b) and (c) below.

If $\ga_1$, $\ga_2$ are flat $\su$ connections over $Y$ and $\ga_1$ is
irreducible then by taking the holonomy of connections along the path
$\R_-\times\{y_0\}$ one obtains a trivialization
\[f_-:M(\ga_1,\ga_2)\times\co^3\oset\approx\to\e(\ga_1,\ga_2)\otimes\co\]
of the natural complex $3$-plane bundle over $M(\ga_1,\ga_2)$.
Similarly, if $\ga_2$ is irreducible one gets a trivialization $f_+$
of $\e(\ga_1,\ga_2)\otimes\co$ in terms of holonomy along $\R_+\times\{y_0\}$.

Now fix linearly independent elements $e_1,e_2\in\co^3$.
If $D_1\subset M(\theta,\al)$, $D_2\subset M(\beta,\theta)$ are compact
sets then by modifying the sections $s_1$, $s_2$ in a small neighbourhood
of $r(D_1\cup D_2)$ (this will not affect the chain map $u$) one can arrange
that
\[s_j(r(A))=\be{cases}
	f_-(A,e_j) & \text{if $[A]\in D_1$}\\
	f_+(A,e_j) & \text{if $[A]\in D_2$.}\end{cases}\]
Here we are making use of the following four facts:
\be{description}
\item[$\bullet$]If $K\subset\cb^*((-1,1)\times Y)$ is the union of all images
$r(M(\ga_1,\ga_2))$ where the flat connections $\ga_j$ are irreducible and
$\dim\,M(\ga_1,\ga_2)\le4$ then $K$ is compact.
\item[$\bullet$]Unique continuation: If two anti-self-dual connections over
$\ry$ are gauge
equivalent over $\band$ then they must be gauge equivalent over $\ry$.
\item[$\bullet$]Restriction to $\band$
defines smooth embeddings of $M(\theta,\al)$
and $M(\beta,\theta)$ into $\cbb$.
\item[$\bullet$]$\cbb$ admits smooth partitions
of unity (recall that we are working
with $L^p_1$ connections with $p>4$ an even integer).
\end{description}

{\bf(II)} We will now state a gluing theorem which describes the elements of 
$M(\beta,\al)$ that are close to the trivial connection over $\band$.

Fix a small positive constant $\eps_1$ and let $U$ be the set of all elements
of $M(\beta,\al)$ which over the band $\band$ can be represented by a 
connection form $a$ with $||a||_{L^2_1}<\eps_1$. As $\eps_1$ becomes
smaller, elements of $U$ will more and more resemble broken gradient lines
from $\al$ to $\beta$ factoring via the trivial connection $\theta$. Hence
for sufficiently small $\eps_1$ there is a natural map
\[U\oset\psi\to\cem(\theta,\al)\times\cem(\beta,\theta).\]

Let $\eps_2$ be another small positive constant and choose a smooth function
$\phi:\R\to\R$ satisfying $\phi'\ge0$, $\phi(t)=1$ for $t\le-1$,
and $\phi(t)=0$ for $t\ge1$.
For any flat $\su$ connection $\ga$ over $Y$ define smooth, real functions
$\tau_1$ on $M(\ga,\al)$ and $\tau_2$ on $M(\beta,\ga)$ implicitly by
\be{align*}
&\int_{\ry}|F(A)_{(t,y)}|^2\phi(\tau_1(A)+t)\,dt\,dy=\eps_2\\
&\int_{\ry}|F(A)_{(t,y)}|^2\phi(\tau_2(A)-t)\,dt\,dy=\eps_2.
\end{align*}

For any real number $T$ let $U_T$ be the subset of $U$ defined by
the inequalities $\tau_j>T$, $j=1,2$. 

For any $[A]\in M(\beta,\al)$ let $\eta(A)\in\so$ be the holonomy of $A$
along the path $\R\times\{y_0\}$ in the positive direction.

In gluing theory (see \cite{DK}, \cite{D5}) one proves the following theorem.

\be{thm}\label{gl-thm}
For sufficiently large $T>0$ (depending on $\eps_1$ and $\eps_2$), the map
\[(\tau_1,\tau_2,\eta,\psi):U_T\to(T,\infty)\times(T,\infty)\times\so
\times\cem(\theta,\al)\times\cem(\beta,\theta)\]
is an orientation preserving diffeomorphism.
\end{thm}

This theorem expresses our convention for relating the
orientations of the moduli spaces $M(\theta,\al)$ to those of the spaces
$M(\beta,\theta)$.

We now continue the proof of Theorem~\ref{chain-map}.

{\bf(III)} This part of the proof is related to the computation of
the coefficient of $\del\otimes\del'$ in the Theorem. Choose smooth maps
\[\zeta_j:\R\times\so\to\co^3,\quad j=1,2\]
satisfying
\[\zeta_j(t,g)=\be{cases}
\bar e_j, & t\le-1\\
g\inv\bar e_j, & t\ge1,
\end{cases}\]
where $\bar e_j\in\co^3$ is a vector close to $e_j$ which will be specified
later.
We may arrange that $\zeta_1$ has no zeros and that if $\uline{\co}^3$
denotes the trivial complex $3$-plane bundle over $\R\times\so$ then
the section $\zeta=\zeta_2\mod\zeta_1$ of the quotient bundle
$\uline{\co}^3/\co\zeta_1$ is transverse to the zero-section.

We will now compute the number of zeros of $\zeta$, counted with sign.
Let $\Si$ be the suspension of $\so$, which is the union of two cones:
$\Si=C_+\cup C_-$.
Let $\Si_0\to\Si$ be the principal $\so$ bundle whose "clutching map"
$C_+\cap C_-=\so\to\so$ is the identity map. Then
\[\#\zeta\inv(0)=-\langle p_1(\Si_0),[\Si]\rangle=-2.\]

{\bf(IV)} Fix $T>0$ such that the conclusion of Theorem~\ref{gl-thm} holds.
Choose a smooth function $w:\R\to\R$ such that
$w(t)=1$ for $t\le T$ and $w(t)=0$ for $t\ge 2T$. Define two real functions
$\tau,\rho$ on $M(\beta,\al)$ by
\[\tau=(\tau_1\inv+\tau_2\inv)\inv,\]
(this will serve as a smooth approximation to $\min(\tau_1,\tau_2)$), and
\[\rho=
\be{cases}
w\circ\tau&\text{on $U$}\\
1&\text{on $M(\beta,\al)\setmin U$}.
\end{cases}\]
We can ensure that $\rho$ is smooth by choosing $T$ so large that
$\tau<T$ on $\prtl U$.
For $j=1,2$ define two sections $\xi_j$ and $\ti\si_j$ of
$\e(\beta,\al)\otimes\co$ by
\be{align*}
\xi_j(A)&=f_-(A,\zeta_j(\tau_1(A)-\tau_2(A),\eta(A)))\\
\ti\si_j&=\rho\si_j+(1-\rho)\xi_j.
\end{align*}
For a generic choice of the $\bar e_j$ and $\zeta_j$'s the sections 
$\ti\si_1,\ti\si_2$ will satisfy the same transversality assumptions as
$\si_1,\si_2$. Thus if $\ti\si$ is defined as $\si$
with $\ti\si_j$ in place of $\si_j$ then $Z=\ti\si\inv(0)$ is an oriented,
smooth, $1$-dimensional submanifold of $M(\beta,\al)$, and $Z$ can be
described as the locus in $M(\beta,\al)$ where
$\ti\si_1$ and $\ti\si_2$ are linearly dependent.

We will now determine the ends of $Z$.
Suppose $[A_n]$ is a sequence in $Z$ which has no
subsequence which converges in $M(\beta,\al)$. After passing to a subsequence
and applying suitable gauge transformations we may assume $A_n$ converges
over compact subsets of $\ry$ to some anti-self-dual connection $A$.
There are now four possibilities:
\be{description}
\item[(a)]The flat limits of $A$ are both irreducible. Then $[A]$ must lie
in a $4$-dimensional moduli space $M(\ga,\al)$ (where $i(\ga)=0$) or
$M(\beta,\ga)$ (where $i(\ga)=5$). Gluing theory tells us that
the corresponding number of ends of $Z$ is 
$\langle(dv-vd)\al,\beta\rangle$.
\item[(b)]$[A]\in M(\theta,\al)$. This means that $\tau_1(A_n)$ stays 
bounded while $\tau_2(A_n)\to\infty$ as $n\to\infty$. Hence
$\xi_j(A_n)\to f_-(A,\bar e_j)$. Moreover, $[A]$ must satisfy
$T\le\tau_1(A)\le2T$. Now the inequalities $T\le\tau_1\le2T$ define a 
compact subset $D_1\subset M(\theta,\al)$, so by (I) we may assume
the sections $s_1,s_2$ are 
chosen such that $\si_j(A_n)\to f_-(A,e_j)$. Together this implies
that $\ti\si_1(A_n)$ and $\ti\si_2(A_n)$ must be linearly independent for
large $n$, contradicting $\ti\si(A_n)=0$. Thus $[A]$ cannot lie in
$M(\theta,\al)$.
\item[(c)]$A\in M(\beta,\theta)$. This is ruled out just like case (b).
\item[(d)]$A$ is trivial. Then $\tau(A_n)\to\infty$, so $\ti\si_j(A_n)=
\xi_j(A_n)$ for large $n$. The corresponding number of ends of 
$Z$ is $-\#\zeta\inv(0)\cdot\langle(\del\al)\del',\beta\rangle$, and by
(III) we have $\#\zeta\inv(0)=-2$.
\end{description}
This completes the proof of Theorem~\ref{chain-map}.\square

\subsection{Cobordisms}

We will now study how the $u$-map is related to maps between Floer groups
induced by cobordisms. 
There are many different cases here that one could 
consider (ie different kinds of cobordisms and bundles etc), and
we will focus on what is in a sense the most general
case that we will encounter.

Let $W$ be a connected Riemannian $4$-manifold with two cylindrical ends,
$\R_-\times Y_1$ and $\R_+\times Y_2$. Suppose $Y_1$ and $Y_2$ are
homology spheres, $H_1(W;\z)=0$, and the intersection form of $W$ is negative
definite. If $\al_1,\al_2$ are flat $\su$ connections over $Y_1,Y_2$,
respectively, not both trivial, let $M(W;\al_2,\al_1)$ denote
the moduli space of anti-self-dual $\su$ connections over $W$ with flat limits
$\al_1$ at $-\infty$ and $\al_2$ at $\infty$, and with dimension in the range
$[0,7]$. There is a degree preserving cochain homomorphism
\[W^*:\cf^*(Y_1)\to \cf^*(Y_2),\quad
\al\mapsto\sum_\beta(\#\,M(W;\beta,\al))\,\beta\]
where the sum is taken over all gauge equivalence classes $\beta$ of flat $\su$
connections over $Y_2$ of the same index as $\al$. 
There is also a homomorphism $\del_W:\cf^5(Y_1)\to\z$ and an element
$\del'_W\in \cf^0(Y_2)$ obtained by counting points in zero-dimensional moduli
spaces over $W$ with trivial limit over one end.

For $j=1,2$ choose a generic pair of sections of the natural complex
$3$-plane bundle $\f_j\otimes\co\to\cb^*(Y_j\times(-1,1))$.
As above this defines
a homomorphism $v:\cf^*(Y_j)\to \cf^{*+4}(Y_j)$.

\be{thm}\label{ucob}
There exists a graded homomorphism $\phi:\cf^*(Y_1)\to \cf^{*+3}(Y_2)$
such that
\[vW^*-W^*v+2(\del_W\otimes\del'+\del\otimes\del'_W)=d\phi+\phi d\]
as homomorphisms $\cf^*(Y_1)\to \cf^{*+4}(Y_2)$,
where $\del:\cf^4(Y_1)\to\z$ and $\del'\in \cf^1(Y_2)$ are as
defined in Section~\ref{Floer-homology}.
\end{thm}
\proof The proof will be quite similar to the proof of Theorem~\ref{chain-map},
so we only indicate the new features. We will use the technique of
``moving the base-point''. In a sense, this will make up for the lack of
translation in moduli spaces over $W$.
To this end, choose base-points
$y_j\in Y_j$, $j=1,2$, and a smooth path $\ga:\R\to W$ such that
\[\ga(t)=\be{cases}
(t,y_1)& \text{for $t\le-1$}\\
(t,y_2)& \text{for $t\ge1$.}
\end{cases}\]
Set $W_0=W\setmin((\infty,-2]\times Y_1\cup[2,\infty)\times Y_2)$.
We may assume that $\ga(t)\in W_0$ for $|t|\le2$. Set $\ga_0=\ga|_{(-2,2)}$.

Now let $\al,\beta$ be flat $\su$ connections
over $Y_1$ and $Y_2$, respectively, not both trivial, and set
$M=M(W;\beta,\al)$. Let $\ub\to M\times W$ and $\ub_0\to
\cb^*(W_0)\times W_0$ be the universal Euclidean $3$-plane bundles
(see \cite{DK}), and let
\[\e=(\text{id}_W\times\ga)^*(\ub),\quad
\e_0=(\text{id}_{W_0}\times\ga_0)^*(\ub_0)\]
be the pull-back bundles over $M\times\R$ and $\cb^*(W_0)\times(-2,2)$,
respectively. Choose a generic pair of sections of $\e_0\otimes\co$.
We can pull back these sections
and the sections of $\f_j\otimes\co$ by the restriction maps
\be{align*}
r_0:M\times(-2,2)&\to\cb^*(W_0)\times(-2,2),\quad([A],t)\mapsto([A|_{W_0}],t)\\
r_1:M\times(-\infty,-1)&\to\cb^*((-1,1)\times Y_1),\quad
([A],t)\mapsto[A|_{(t-1,t+1)\times Y_1}]\\
r_2:M\times(1,\infty)&\to\cb^*((-1,1)\times Y_2),\quad
([A],t)\mapsto[A|_{(t-1,t+1)\times Y_2}]
\end{align*}
to obtain pairs of sections of $\e\otimes\co$ over $M\times(-2,2)$,
$M\times(-\infty,-1)$, and $M\times(1,\infty)$, respectively.
Piecing these together
using a partition of unity we obtain two sections $\si_1,\si_2$ of 
$\e\otimes\co$. If $\dim\,M\le4$ then we may assume $\si_1$ has no zeros
and that the section $\si=\si_2\mod\si_1$ of the quotient bundle
$(\e\otimes\co)/\co\si_1$ is transverse to the zero-section. 
If $\al,\beta$ are both irreducible and
$\dim\,M=3$ then $\si\inv(0)$
is a finite set of oriented points, and we define the matrix coefficient
$\langle\phi(\al),\beta\rangle$ by
\[\langle\phi(\al),\beta\rangle=\#\si\inv(0).\]

The remainder of the proof follows the proof of Theorem~\ref{chain-map}
quite closely. Let $\al$ and $\beta$ be irreducible and
$i(\beta)\equiv i(\al)+4\mod8$, so that $\dim\,M=4$.
We first note that, as in (I), certain alterations may be made to the sections
of $\f_j\otimes\co$ used above, without affecting the definition of the
$u$-map on $\cf^*(Y_j)$.
Also, there is an analogue of Theorem~\ref{gl-thm} for the present setup
which we use to redefine the sections 
$\si_1,\si_2$ for elements $([A],t)$ of $M\times\R$
where either $t<<0$ and $r_1(A)$ is close to the trivial connection, or
$t>>0$ and $r_2(A)$ is close to the trivial connection.

This being done, let $Z\subset M\times\R$ be the zero-set of
$\ti\si$ (the modification of $\si$). Thus $Z$
is an oriented, smooth, $1$-dimensional submanifold of $M\times\R$.
To describe the ends of $Z$, suppose $([A_n],t_n)$ is a sequence in $Z$
which has no convergent subsequence in $M\times\R$. After passing to a
subsequence we may assume $t_n$ has a limit $L$ in $[-\infty,\infty]$.

If $L\in\R$ then we may pass to a subsequence in which $A_n$ converges 
modulo gauge transformations over compact subsets of $W$ to some
$[A]\in M(W;\ga_1,\ga_2)$. For dimensional reasons $[A]$ cannot be reducible.
Hence for transversality reasons we must have
$\dim\,M(W;\ga_1,\ga_2)=3$, so either $\ga_1=\al$ or $\ga_2=\beta$.
The corresponding number of ends of $Z$ is $-\la(d\phi+\phi d)\al,\beta\ra$.

The case $L=-\infty$ is analogous to the proof of Theorem~\ref{chain-map}, and
the corresponding number of ends is
$\la(-W^*v+2\del\otimes\del'_W)\al,\beta\ra$.
Similarly, in the case $L=\infty$ the number of ends is
$\la(vW^*+2\del_W\otimes\del')\al,\beta\ra$.\square

\subsection{Reduced Floer groups}

Let $Y$ be an oriented homology $3$-sphere. If $n$ is a non-negative integer 
then $v^n\del'$ lies in $\cf^{1+4n}(Y)$, hence by Theorem~\ref{chain-map},
\[dv^n\del'=v^nd\del'=0,\]
and similarly $\del v^nd=0$. Let
\[\del'_n\in \hf^{1+4n}(Y)\]
be the cohomology class of $v^n\del'$, and let
\[\del_n:\hf^{4-4n}(Y)\to\z\]
be the homomorphism induced by $\del v^n$.

It follows from the cochain homotopy formula in Theorem~\ref{chain-map} (ii)
that either $\del_0$ or $\del'_0$ must be zero. Moreover, if $\del_0$
is zero then $\del_n$ vanishes for all $n$, and similarly if $\del'_0$ is zero.

In general, we do not expect that $\del_n$ is a topological
invariant of $Y$ for $n>1$. However, the following theorem shows 
that $\del_n$ is a topological invariant 
modulo $\del_0,\dots,\del_{n-1}$.
Recall that one can compare Floer groups defined by 
different metrics and perturbations by means of the homomorphism
induced by a cobordism $\ry$ where the metric and perturbation 
interpolates between the given ones on $Y$ (see \cite{F1}).
Thus the question is how $\del_n$ behaves under maps induced
by cobordisms. 

\be{lemma}\label{del-funct}If $W,Y_1,Y_2$ are as in Theorem \ref{ucob} then
\[\del W^*=\del+\del_Wd\]
as maps $\cf^4(Y_1)\to\z$.
\end{lemma}
\proof Let $\al$ be an irreducible flat 
$\su$ connection of index~$4$ over $Y_1$. We will determine the ends
of the $1$-dimensional $\su$ moduli space $M=M(W;\theta,\al)$. Let $[A_n]$ be
a sequence in $M$.  
By taking a subsequence we may arrange that $[A_n]$ converges modulo gauge
transformations to some instanton $A$ over $W$. For index reasons $A$ must
be either irreducible or trivial.

If $A$ is irreducible then it must have
index $0$, and factorization has occured through an irreducible
flat connection of index $4$ over $Y_1$ or $Y_2$.
The corresponding number of ends of $M$ is $(\del W^*-\del_Wd)\al$.

The number of ends of $M$ corresponding to the case when $A$ is trivial 
is $-\del\al$. Here we are using the assumption that $b^+_2(W)=0$, which
just means that the trivial connection over $W$ is a regular solution of the
instanton equation.\square

\be{thm}\label{del-cob}
Let $W,Y_1,Y_2$ be as in Theorem \ref{ucob}. Then $\del_0W^*=\del_0$
and $W^*\del'_0=\del'_0$. More generally, there are integers $a_{ij},b_{ij}$
such that for any non-negative integer $n$,
\[\del_nW^*=\del_n+\sum_{i=0}^{n-1}a_{in}\del_i,\qquad
W^*\del'_n=\del'_n+\sum_{i=0}^{n-1}b_{in}\del'_i.\]
\end{thm}

\proof The statement $\del_0W^*=\del_0$ follows from Lemma~\ref{del-funct},
and the proof that $W^*\del'_0=\del'_0$ is similar.
The remainder of the theorem is then
a simple consequence of Theorem~\ref{ucob}.\square

Of course, one can take $a_{in}=0=b_{in}$ when $i$ and $n$ have different
parity. 

Let $B^*\subset \hf^*(Y)$ be the linear span of the classes $\del'_n$, $n\ge0$.
Thus $B^1$ is spanned by the $\del'_{2n}$, and $B^5$ is spanned by
the $\del'_{2n+1}$. If $q\not\equiv1\mod4$ then $B^q=0$. Also, set
\[Z^*=\bigcap_{n\ge0}\ker(\del_n)\subset \hf^*(Y),\]
where by definition $\del_n$ is zero in degrees different from
$(4-4n)$ mod~$8$.
\be{defn}The reduced Floer group $\hfh^*(Y)$ is defined by
\[\hfh^q(Y)=Z^q/B^q.\]
\end{defn}
We will now define the $u$-map on the reduced Floer group. Let
$\ker(d)\subset \cf^*(Y)$ be the Floer cocycles, and let
$\pi:\ker(d)\to \hf^*(Y)$. Using
Theorem~\ref{chain-map} it is easy to check that $v$ maps $\pi\inv(Z^*)$
into itself and $\pi\inv(B^*)$ into itself, hence $v$ induces a degree
$4$ endomorphism $u$ of 
\[\hfh^*(Y)=\pi\inv(Z^*)\,/\,\pi\inv(B^*).\]
\be{thm}\label{ucob2}
Let $W,Y_1,Y_2$ be as in Theorem~\ref{ucob}. Then $(W^*)\inv(Z^q)=Z^q$
and $W^*(B^q)=B^q$ for every $q$. In particular, $W$ induces a homomorphism
$\hfh^*(Y_1)\to\hfh^*(Y_2)$. Moreover, this homomorphism
commutes with the $u$-maps.
\end{thm}
\proof This follows immediately from Theorems~\ref{ucob} and \ref{del-cob}.
\square

Just as for the ordinary Floer groups the map between the reduced Floer groups
induced by $W$ is independent of the metric on $W$, as long as the metric is
on product form on the ends.

\be{cor}The $\z[u]$ module $\hfh^*(Y)$ is an invariant of the
oriented, smooth manifold $Y$.
\end{cor}

\section{Reducible connections}\label{red-conn}

We will now see how one can obtain information about the homomorphisms
$\del_n:\hf^{4-4n}(Y;\q)\to\q$ when the oriented homology sphere $Y$ bounds a 
smooth, compact, oriented $4$-manifold with non-standard definite intersection
form. The main ideas here are due to Donaldson~\cite{D5}. It should be stressed
that our results do by no means give a complete description, since we only take
into account the lowest stratum of reducibles. 

Let $X$ be a smooth, oriented Riemannian $4$-manifold with
one cylindrical end $\rpy$, where $Y$ is an integral
homology sphere. Let $E\to X$ be a  $U(2)$ bundle and fix a smooth connection
$A^{\det}$ in $\La^2(E)$ which is trivial over the end. For any integer $k$ let
$M(E,k)$ denote the moduli space of (projectively) anti-self-dual connections
$A$ in $E$, with central part $A^{\det}$, which are asymptotically trivial
over the end and satisfy
\[\fact\int_X\text{tr}(F(A)^2)=k-\frac12c_1(E)^2.\]
Here $F(A)$ is the curvature of $A$. Then the expected 
dimension of $M(E,k)$ is
\[\dim M(E,k)=8k-2c_1(E)^2-3(1-b_1(X)+b_2^+(X)).\]
Note that if $F_0(A)$ is the traceless part of $F(A)$ then
\be{equation*}
0\le\fact\int_X\text{tr}(F_0(A)^2)=k-\frac14c_1(E)^2.
\end{equation*}
The inequality follows because $\text{tr}(F_0(A)^2)=|F_0(A)|^2$ when
$F_0^+(A)=0$.

\be{lemma}\label{red_class}
If $X$ is connected and $b_1(X)=0$ then the set of
reducible points in the moduli
space $M=M(E,k)$ is in one-to-one correspondence
with the set of unordered pairs
\be{equation*}
R=\{\{z_1,z_2\}\subset H^2(X;\z)\st(z_1-z_2)_{\R}\in\harm^-;\;
z_1+z_2=c_1(E);\;z_1z_2=k\}.
\end{equation*}
Here $(\cdot)_{\R}$ denotes the real reduction of an integral cohomology class
and $\harm^-$ is the space of anti-self-dual closed $L^2$ $2$-forms
on $X$.
\end{lemma}

\proof We define a map $r:M^{\text{red}}\to R$, where $M^{\text{red}}$
is the set of reducible points in $M$. If $u$ is any
automorphism of $E$ such that $u\not\in U(1)$ and $u(A)=A$ then $A$ preserves
the eigenspace decomposition $E=L_1\oplus L_2$. This splitting is unique 
up to order unless $A$ is projectively trivial, in which case all
$A$-invariant rank~$1$ sub-bundles of $E$ are isomorphic.

Set $z_j=c_1(L_j)$ and $r(A)=\{z_1,z_2\}$. Then $r$ is well defined. For
instance, to verify $z_1z_2=k$ let $F(A|_{L_j})=i\phi_j\in i\Omega^2(X;\R)$,
which represents $2\pi ic_1(L_j)$. Since the
anti-self-dual closed $L^2$ forms $\phi_j$ decay exponentially on the end, we
have
\[k-\frac12(z_1+z_2)^2=\fact\int_X\text{tr}(F(A)^2)=\fact\int_X(-\phi_1^2
-\phi_2^2)=-\frac12(z_1^2+z_2^2),\]
hence $k=z_1z_2$ as claimed. 

To see that $r$ is a bijection recall that by Hodge theory a
$U(1)$ bundle $L\to X$ admits a finite energy anti-self-dual connection
precisely when the $L^2$-harmonic form representing
$c_1(L)$ is anti-self-dual, and in that case the connection is unique up
to gauge since $b_1(X)=0$.
\square

With $X$ as above, let $\ct\subset H^2(X;\z)$ be the torsion subgroup and set
$\cl=H^2(X;\z)/\ct$. Let $\ti R$ be the set of
{\em ordered} pairs $(z_1,z_2)$ such that
$\{z_1,z_2\}\in R$. Then the torsion subgroup $\ct$ of $H^2(X;\z)$
acts freely on $\ti R$
by $t\cdot(z_1,z_2)=(z_1+t,z_2-t)$ for $t\in\ct$.
By associating to a pair $(z_1,z_2)$ the class $(z_1-z_2))_{\R}$ we obtain
a natural identification
\be{equation}\label{alt:r}
\ti R/\ct=\{z\in(c+2\cl)\cap\harm^-(g)\st z^2+4k=c^2\},
\end{equation}
where $c=c_1(E)$.
Note that if $c\not\in2H^2(X;\z)$ then
the projection $\ti R\to R$ is two-to-one.

By a {\em lattice} we shall mean a finitely generated free abelian group
$\cl$ with a non-degenerate symmetric bilinear form $b:\cl\times \cl\to\z$.
For $x,y\in \cl$ we usually write $x\cdot y$ for the pairing
$b(x,y)$, and $x^2$ instead of $x\cdot x$. The dual lattice
$\text{Hom}(\cl,\z)$ will be denoted $\cl^\#$.

\be{defn}\label{eta_def}
Let $\cl$ be a (positive or negative) definite lattice. A vector
$w\in \cl$ is called {\em extremal} if $|w^2|\le|z^2|$ for all $z\in w+2\cl$.
If $w\in \cl$, $a\in \cl^\#$, and $m$ is a non-negative integer
satisfying $w^2\equiv m\mod2$ set
\[\eta(\cl,w,a,m)=\sum_{\pm z}(-1)^{((z+w)/2)^2}(a\cdot z)^m,\]
where the sum is taken over all unordered pairs $\{z,-z\}\subset w+2\cl$
such that $z^2=w^2$.
If $m=0$ then we write $\eta(\cl,w)=\eta(\cl,w,a,m)$.
\end{defn}

\be{prop}\label{red_count}
Let $X$ be a smooth, compact, oriented $4$-manifold
with a homology sphere $Y$ as boundary and with $b_1(X)=0$.
Suppose the intersection form on
$\cl=H^2(X;\z)/\ct$
is negative definite,  where $\ct$ is the torsion subgroup.
Let $c\in H^2(X;\z)$, $a\in H_2(X;\z)$, and
let $m$ be a non-negative integer such that
$c^2\equiv m\mod2$ and $-c^2\ge2$.
Set $n=-(c^2+m)/2-1$. If $c$ reduces to a non-zero 
extremal vector $w\in\cl$ then the Donaldson invariant
$D_X^c(a^m)\in\hf^{4-4n}(Y;\q)$ is well defined and
\[\del_0u^j\cdot D_X^c(a^m)=\be{cases}
0&\text{for $0\le j<n$}\\
\pm2^{-m}\,|\ct|\,\eta(\cl,w,a,m)&\text{for $j=n$.}
\end{cases}\]
\end{prop}
Here we set $\del_0x=0$ if $x\in\hf^0$.

\be{cor}[Donaldson]\label{ex-int}
If the intersection form of $X$ is not diagonal over the integers
then $\del_0:\hf^4(Y;\q)\to\q$ is non-zero.
\end{cor}
\proof Let $\bar\cl\subset\cl$ be the orthogonal complement of 
all vectors of square $-1$. The assumption is that
$\bar\cl\neq0$. 
Let $w\in\bar\cl$ be a non-zero vector of smallest length,
and set $m=-w^2-2$. Choose a class $a\in H_2(X;\z)$ with $a\cdot w=1$.
Then $\eta(\cl,w,a,m)=1$, so $\del_0\neq0$ by the proposition.\square

{\em Proof of Proposition.}
First add a half-infinite cylinder $\rpy$ to $X$ and
choose a cylindrical end
metric on this new manifold (also denoted $X$). 
Let $E\to X$ be the $U(2)$ bundle with $c_1(E)=c$, and set
$M_k=M(E,k)$.

By Lemma~\ref{red_count}, the moduli space $M_k$ contains no reducibles
if $k<0$. Now consider $k=0$. Let $R$ and $\ti R$ be
as above. It is convenient to fix an ordering of each
pair $\{z_1,z_2\}\in R$.
By making a small perturbation of the anti-self-dual equations
near each reducible point as in 
\cite{D1} we can arrange that for each
$\bz=(z_1,z_2)$ the corresponding reducible point in $M_0$ has
an open neighbourhood $C_\bz$ in $M_0$
which is homeomorphic to a cone on some complex projective space $P_\bz$.
Note that reversing the order of $z_1,z_2$
reverses the complex structure on $P_\bz$.

The boundary orientation that
$P_\bz$ inherits from $M_0^\#=M_0\setmin\cup_\bz C_\bz$ differs from
the complex orientation of $P_\bz$ by a sign $\eps(z_1,z_2)$. If 
$(\bar z_1,\bar z_2)$ is another element of $\ti R$ then it follows from
\cite{D2} that
\[\eps(z_1,z_2)\eps(\bar z_1,\bar z_2)=(-1)^{(z_1-\bar z_1)^2}.\]
Now consider the universal $\so$ bundle 
\[\e\to P_\bz\times X\]
and define the $\mu$-map $H_i(X)\to H^{4-i}(P_\bz)$ as usual by
$\mu(b)=-\frac14p_1(\e)/b$.
If $e\in H^2(P_\bz)$ is the Chern class of the tautological 
line bundle then we have
\be{equation*}
\mu(1)=-\frac14e^2;\qquad\mu(a)=-\frac12\langle z_1-z_2,a\rangle e
\end{equation*}
for any $a\in H_2(X;\z)$ (see \cite{DK}).

Let $M_0'$ be the oriented $1$-manifold with boundary obtained by
cutting down $M_0^\#$
according to the monomial $x^ja^m\in\mathbb A(X)$ as in Section~\ref{Floer-homology}.
The boundary points of $M_0'$ lie on the links of the reducibls points, while
the ends correspond to factorizations on $\rpy$ through flat connections
of index $4$. Counted with sign, the number of boundary points in $M_0'$ plus
the number of ends must be zero. This gives
\[\del D_X^c(x^na^m)=\pm2^{-2n-m}\,|\ct|\,\eta(\cl,w,a,m),\]
where the left hand side is the number of ends. The invariant
$ D_X^c(x^na^m)\in \hf^4(Y)$ is well defined because $M_k$ contains
no reducibles when $k<0$. The same argument shows that 
$\del D_X^c(x^ja^m)=0$ for $0\le j<n$.
The Proposition now follows because 
\[ D_X^c(x^ja^m)=(\frac14u)^j D_X^c(a^m)\qquad\text{for $0\le j\le n$}.\]
One can prove this by ``moving the base-point''
along the path $[0,\infty)\times Y$ (as in Section~\ref{umap}).
However, to run this argument one needs to know that no irreducible
connections in $M_k$ can restrict to
a reducible connection over the end $\rpy$. In other words, $M_k$ must not
contain any twisted reducibles (cfr Lemma~4.3.21 in \cite{DK}).
Fortunately, this holds at
least generically:

After perturbing the metric on $X$ in a small ball we may
assume that there are no non-flat twisted reducibles
in $M_k$ for any $k$ (see \cite{KM3}, Section 2(i)).
The main point here is that if $\lla$ is a real line bundle over $X$ and
$b_i(\lla)$ are the Betti numbers of $X$ with coefficients in $\lla$ then
as pointed out in \cite{KM3},
\[b_0(\lla)-b_1(\lla)+b_2^+(\lla)\]
is the same for all line bundles $\lla$. If $b_1(X)=0$ and $\lla$ is 
non-trivial this gives $b_2^+(\lla)>b_2^+(X)=0$.

Moreover, since $w$ is not divisible
by $2$ in $\cl$, there are no flat reducibles in $M_k$. Any flat, twisted
reducible connection in $M_k$ must therefore be irreducible. But
$M_k$ can only contain a flat connection for $k=c_1(E)^2/4$; in that case
$M_k$ has expected dimension $-3$ and is generically empty
(compare \cite{D2}).
Hence we may assume there are no twisted reducibles in $M_k$ for any $k$.
In particular, no irreducible connection in $M_k$ will restrict to a reducible
connection over the end. \square

We will now apply this proposition to obtain information about the 
Floer cohomology of the \poinc\ sphere $S=\Si(2,3,5)$.
It is well known that the
Milnor fibre of the $E_8$ singularity can be smoothly embedded in a $K3$ 
surface $X$. This gives a splitting $X=X_1\cup_S X_2$, where 
the intersection form of $X_1$ is $-E_8$ and the intersection form of
$X_2$ is $-E_8\oplus3
\bigl(\be{smallmatrix}
0&1\\1&0\end{smallmatrix}\bigr)$. If $e_1,\dots,e_8$ is an
orthonormal basis for $\R^8$ then
\be{equation*}
E_8=\{\textstyle{\sum\limits_{i=1}^8} x_i e_i\:|\: 2x_i\in\z;\;
x_i-x_j\in\z;\; \textstyle{\sum\limits_{i=1}^8}x_i\equiv0\,(2)\}.
\end{equation*}
For $j=1,2$, let $w_j\in H^2(X_j;\z)$ be the element corresponding to
$e_1+e_2\in E_8$, and let $z_j\in H^2(X_j;\z)$ be the element 
corresponding to $e_1+e_2+e_3+e_4\in E_8$. Recall that the relative invariants
of $X_1$ and $X_2$ take values in the Floer cohomology and homology of $S$,
respectively. By gluing theory and knowledge of
a certain Donaldson invariant of $X$ (see \cite{Kron1},
or \cite{DK}, Proposition~9.1.3) we have
\be{align*}
 D^{w_1}_{X_1}(1)\cdot D^{z_2}_{X_2}(1)&=D^{w_1+z_2}_X(1)=\pm1;\\
 D^{z_1}_{X_1}(1)\cdot D^{w_2}_{X_2}(1)&=D^{z_1+w_2}_X(1)=\pm1,
\end{align*}
where the signs depend on the homology orientations of $X_1$ and $X_2$.

But $S$ has precisely two equivalence classes of irreducible flat $\su$
connections, which are both non-degenerate, so $\al= D^{w_1}_{X_1}(1)$ and 
and $\beta= D^{z_1}_{X_1}(1)$ must be generators of
$\hf^4(S;\z)=\z$ and $\hf^0(S;\z)=\z$, respectively.

It is easy to check that $e_1+e_2$ and $e_1+e_2+e_3+e_4$ are extremal vectors
in $E_8$ satisfying
\be{gather*}
\eta(E_8,e_1+e_2)=1;\\
\eta(E_8,e_1+e_2+e_3+e_4)=8.
\end{gather*}
Proposition~\ref{red_count} (with $m=0$ and $a=0$) now gives 
$\del_0\al=\pm1$ and $\del_0u\beta=\pm8$, and we deduce
the following proposition.
\be{prop}\label{poinc_cup}
For the \poinc\ sphere $S=\Si(2,3,5)$ the following holds.
\be{description}
\item[(i)]$\del_0:\hf^4(S;\z)\to\z$ is an isomorphism.
\item[(ii)]$u:\hf^0(S;\z)\to \hf^4(S;\z)$ is multiplication by $\pm8$.
\end{description}
\end{prop}
This result was first proved by Kronheimer~(unpublished) and
Austin~\cite{Au1}.

It follows from the proposition that the reduced Floer group of the 
\poinc\ sphere is zero.

\section{Floer's exact triangle}\label{exact_tri}

This section gives a brief description of Floer's exact triangle;
for further details we refer to the exposition \cite{BD1}.
Let $Y_0$ be a closed, oriented $3$-manifold and
$P_0\to Y_0$ an admissible principal $\so$ bundle.
Given any orientation
preserving $\so$ equivariant embedding
\be{equation}\label{ka_map}
\ka_0:D^2\times S^1\times\so\to P_0
\end{equation}
we can form the surgery cobordism
\[Q_0=\oline{D^2\times D^2}\times\so\cup_\ka[0,1]\times P_0,\]
where $\ka$ maps into $\{1\}\times P_0$. 
Here the bar means we take the opposite of the standard orientation.
The oriented boundary of $Q_0$ is
$P_1\cup\oline P_0$, where $P_1$ is the result of the
surgery on $P_0$ determined by $\ka_0$. Let $f:S^1\to\so$ be a
homotopically non-trivial map and define an equivariant embedding
\be{align*}
\ka_1:D^2\times S^1\times\so&\oset\approx
\to S^1\times D^2\times\so\subset
P_1,\\
(z,w,u)&\mapsto(w\inv,zw\inv,f(w)u).
\end{align*}
Here we regard $D^2$ as the unit disc in $\co$, and $S^1=\prtl D^2$.
Iterating this process we get a sequence $(P_n,\ka_n)$, 
$n=0,1,2,\dots$. The bundle $P_{n+1}$
is obtained from $P_n$ by cutting out $\im(\ka_n)$ and re-gluing
using a certain equivariant diffeomorphism $\xi$ of $S^1\times S^1
\times\so$. A crucial point here is that $\xi^3$ covers the identity
map on $S^1\times S^1$, and its associated map $S^1\times S^1\to\so$
is null-homotopic. Thus we may identify 
\[(P_{n+3},\ka_{n+3})=(P_n,\ka_n).\]

From now on we assume all bundles $P_n$ are admissible. There are
then two possibilities: either (i) $Y_n$ is not a homology sphere
for any $n$, or (ii) for some $n$, the manifold
$Y_n$ is a homology sphere, while $Y_{n+1}$ and $Y_{n+2}$ are
the result of $-1$ surgery and $0$ surgery (respectively)
on the knot in $Y_n$ determined by $\ka_n$.
In both cases the surgery cobordism $Q_n$ from
$P_n$ to $P_{n+1}$ induces a homomorphism $\al_n:\hf^*(P_n)\to
\hf^*(P_{n+1})$.
Floer's theorem now says that for every $n$, the composite
homomorphism $\al_{n+2}\al_{n+1}\al_n$ shifts degrees by $-3\mod8$,
and
\[\im(\al_n)=\ker(\al_{n+1}).\]
(It is easy to compute the shift in degrees in case (ii) above;
on the other hand
the shift must always be the same, as can be seen from the addition 
property for the index over $4$-manifolds with tubular ends.) If $Y_0$ and
$Y_1$ are homology spheres and we give $Y_2$ the spin structure that extends
over $W_1$ then the long exact sequence takes the form
\[\dots\to \hf^q(Y_0)\to \hf^q(Y_1)\to \hf^q(Y_2)\to \hf^{q-3}(Y_0)\to\dots,\]
where we have written $\hf^q(Y_2)$ instead of $\hf^q(P_2)$.

\section{The nilpotency of $u^2-64$}\label{nilpot}

In this section we establish the nilpotency of $u^2-64$ in both 
the Floer group of a non-trivial, admissible $\so$ bundle and in the 
reduced Floer group of a homology sphere.
In both cases the proof of nilpotency of $u^2-64$
begins by representing the $3$-manifold as surgery on an 
algebraically split framed link in $S^3$. We then use Floer's exact
triangle and a link reduction scheme. A central part in the reduction
argument is that if $F$ is the non-trivial $\so$-bundle over the $2$-torus
then $u^2=8^2$ in $\hf^*(S^1\times F)$. The number $8$ is derived from
Proposition~\ref{poinc_cup}~(ii).

We first discuss non-trivial $\so$ bundles.
\be{prop}\label{fintyp1a}
Let $Y$ be a closed oriented $3$-manifold with $H_1(Y;\z)$ torsion free
and non-zero. Let $P\to Y$ be a non-trivial $\so$ bundle.
Then for some $n\ge1$,
\[(u^2-64)^n=0\]
as an endomorphism of $\hf^*(P)$.
\end{prop}

Before proving this proposition we will deduce from it a stronger result.
Let $\Si$ be a surface of genus $g$ and
$F\to\Si$ the non-trivial
$\so$ bundle. Consider the Floer cohomology $\hf^*_g$ of the $\so$ bundle 
$S^1\times F\to S^1\times\Si$.
Let $\psi\in \hf^*_g$ be the element obtained
by counting points (with signs) in zero-dimensional instanton moduli spaces
in the bundle $D^2\times F\to D^2\times\Si$
(adding a tubular end to $D^2\times\Si$ as usual).
For $g\ge0$ let $N_g$ be the smallest
non-negative integer $n$ such that $(u^2-64)^n\,\psi=0$.
(For rational coefficients it should in principle be
possible to compute all the constants
$N_g$ from \cite{Mun1}. For instance one has $N_1=N_2=1$.)

\be{thm}\label{fintyp1b}
Let $Y$ be a closed, oriented $3$-manifold with
$b_1(Y)>0$, and let $P\to Y$ be an $\so$ bundle such that 
$\R\times P$ is non-trivial over some surface $\Sigma\subset
\R\times Y$ of genus $g$. Then
\[(u^2-64)^{N_g}=0\]
as an endomorphism of $\hf^*(P)$.
\end{thm}

{\em Proof of Theorem \ref{fintyp1b} assuming Proposition~\ref{fintyp1a}:}
For any non-negative integer $n$ let 
\[K_n:\hf^*(P)\to\hf^{*+4n}(P)\]
be the homomorphism defined by cutting down $(4n)$-dimensional moduli spaces
$M(\beta,\al)$ according to the monomial $(4x)^n$ as in
Section~\ref{Floer-homology}, where $x\in H_0(\ry)$ is the point class.
By moving one base-point along a path $[0,\infty)\times\{y_0\}$
as in the proof of Theorem~\ref{ucob} one finds that
$K_n=uK_{n-1}$ for $n\ge1$, so by induction,
\[K_n=u^n.\]

Let $f:D^2\times\Si\to\ry$ be a smooth embedding which maps $\{0\}\times\Si$
onto $\Si$. Let $\{g_t\}_{t\ge0}$ be a 
smooth family of metrics on $\ry$ which stretches $\ry$ along
$f(S^1\times\Si)$. More precisely, let $A\subset D^2$ be an annulus
about the origin and set $U=f(A\times\Si)$. Then
\be{description}
\item[$\bullet$]$g_0$ should be a product metric on $\ry$,
\item[$\bullet$]$g_t$ should be independent of $t$ outside $U$,
\item[$\bullet$]under the identification $U\approx[0,1]\times S^1\times\Si$ the
restriction of $g_t$ to $U$ should have the form $t^2dr^2+ds^2$ for
$t\ge1$, where $r$ is the coordinate on $[0,1]$ and $ds^2$ a fixed 
metric on $S^1\times\Si$.
\end{description}

Now let $K_{n,t}$ be defined as $K_t$ above but using the metric $g_t$,
and such that the geometric representatives for $4\mu(x)$ are obtained by
restricting instantons to fixed subsets of $(\ry)\setmin U$. As is well known,
$K_{n,t}$ is independent of $t$, because the cochain map that defines
$K_{n,t}$ is independent of $t$ up to cochain homotopy. We will now describe
$K_{n,t}$ for large $t$. Let $Q$ denote the restriction of $\R\times P$
to $f(S^1\times\Si)$, and let
$W$ be the manifold $(\ry)\setmin\Si$ with a metric which is on
product form on the end $\R_-\times S^1\times\Si$ and agrees with $g_0$
outside $U\setmin\Si$. Then moduli spaces in $(\R\times P)|_W$ cut down
according to $(4x)^n$ define a homomorphism
\[L_n:\hf^*(Q)\to\text{End}(\hf^*(P))\]
such that for large $t$ one has $K_{n,t}=L_n(\psi)$.

For the final step in the proof we move one base-point along a path
$\R_-\times\{z_0\}$ in $\R_-\times(S^1\times\Si)$ to deduce
\[L_n(\rho)=L_{n-1}(u\rho)\]
for any $\rho\in\hf^*(Q)$ and $n\ge1$. By induction on $n$,
\[L_n(\rho)=L_0(u^n\rho).\]

Putting all this together we obtain
\[(u^2-64)^n=L_0((u^2-64)^n\psi)\]
as an endomorphism of $\hf^*(P)$. Therefore the theorem follows from the
proposition.\square

The proof of Proposition~\ref{fintyp1a} begins with three lemmas.

\be{lemma}\label{repr}
Let $Y$ be a closed oriented $3$-manifold with $H_1(Y;\z)$ torsion free.
Then $Y$ can be represented by a framed link in $S^3$
whose linking matrix is
diagonal. The entries on the diagonal are either
$0$ or $\pm1$, and the number of zeros is $b_1(Y)$.
\end{lemma}
\proof Let $L$ be any framed link in $S^3$ representing $Y$. Thus if $X$ is
the $4$-manifold obtained by attaching $2$-handles to the $4$-ball according
to $L$ then $Y=\prtl X$. Let $j:H_2(Y)\to H_2(X)$ be
the map induced by inclusion. Then the intersection form on
$H_2(X)$ descends to a unimodular form
\[q:H_2(X)/\text{im}(j)\to\z.\]
After adding an unknot with framing $\pm1$ if necessary, we may assume
$q$ is odd and indefinite. By the classification of such forms
\cite{HM} $H_2(X)/\text{im}(j)$ then has
a basis over $\z$ with respect to which $q$ is
diagonal. We can therefore modify $L$ by a sequence of
Kirby moves ${\cal O}_2$ (see \cite{Kirby1})
to obtain a framed link, also representing $Y$,
whose linking matrix is diagonal.\square

\be{lemma}\label{gen1}
Let $P\to Y$ be a non-trivial admissible $\so$ bundle.
If $\R\times P$ is non-trivial over some
embedded torus $T\subset\ry$ then the cup product $u$ on
$\hf^*(P)$ satisfies $u^2=64$.
\end{lemma}

\proof Let $F\to T$ be the non-trivial $\so$ bundle over
the $2$-torus.
Then $\hf^*(S^1\times F)$ has rank $1$ in two degrees differing by $4$, and
is zero in the remaining degrees. Let $\tau$ be the natural 
involution of degree $4$ on $\hf^*(S^1\times F)$. Since $u$ and $\tau$
commute there is a constant $c$ such that
\[u=c\tau.\]
Stretching $\ry$ along $S^1\times T$ we find as above that
$u^2=c^2$ on $\hf^*(P)$. To compute $c$, let
$S$ and $S'$ be the result of $-1$ and $0$ surgery on the negative
$(2,3)$ torus knot in $S^3$, respectively. The exact
triangle provides an isomorphism $\hf^*(S)\oset\approx
\to \hf^*(S')$ which commutes with the cup product
$\hf^0\to \hf^4$. But for $S$ this cup product is multiplication
by $\pm8$, according to Proposition~\ref{poinc_cup}, so $c^2=64$.\square

Before stating the next lemma we observe that if
$Y$ is a closed, oriented $3$-manifold and
$\ga\subset Y$ a null-homologous knot, then any
$\so$ bundle over the complement of
$\ga$ has an (up to isomorphism) unique extension to an $\so$ bundle over $Y$.

\be{lemma}\label{gen1_surg}Suppose Proposition~\ref{fintyp1a}
holds for $P\to Y$. Let $\ga\subset Y$ be a knot which
bounds a surface of genus $1$. If $Y'$ is the result of
$\pm1$ surgery on $\ga$ then
Proposition~\ref{fintyp1a} also holds for
the inherited bundle $P'\to Y'$.
\end{lemma}

{\em Proof.} Since Proposition~\ref{fintyp1a} is insensitive to
the orientation of $Y$ it suffices to consider $+1$ surgery on
$\ga$. Choose an embedding $\beta:D^2\times S^1\oset{\approx}\to N\subset Y$
onto a tubular neighbourhood of $\ga$, and a surface $Z\subset Y$ of
genus $1$, such that
\[\prtl Z=Z\cap N=\beta(\{w\}\times S^1)\]
for some $w\in S^1$. Also, choose a trivialization
of $P|_N$ that does not
extend over $N\cup Z$. This trivialization together with the 
embedding $\beta$ determines a map $\ka_0$ as in \Ref{ka_map}.
In the notation of Section~\ref{exact_tri}, the manifold
$Y_1$ is obtained by $0$ surgery on
$\ga$ and $P_1\to Y_1$ is non-trivial over the torus that one
gets by closing up $Z$ with a disc. Furthermore,
$Y_2$ is the result of $+1$ surgery on $\ga$.
Floer's theorem now provides an exact sequence
\[\hf^*(P_1)\oset{\al_1}\to \hf^*(P_2)\oset{\al_2}\to \hf^*(P_0).\]
Suppose there is a positive integer $n$ such that
$(u^2-64)^n=0$ on $\hf^*(P_0)$. For any $x\in \hf^*(P_2)$ the class
$y=(u^2-64)^nx$ will then lie in the kernel of $\al_2$.
Hence $y=\al_1z$ for some
$z\in \hf^*(P_1)$. By Lemma~\ref{gen1} we have $(u^2-64)z=0$, so
\[(u^2-64)^{n+1}x=\al_1(u^2-64)z=0.\quad\square\]

{\em Proof of Proposition~\ref{fintyp1a}.} We first show that if
the proposition holds when
$b_1(Y)=r-1\ge1$ then it also holds
when $b_1(Y)=r$. Let $Y$ be represented by a framed link $L$ as
in Lemma~\ref{repr} and let $L_1,\dots,L_r$ be the components
with framing $0$, where $r\ge2$. If $\ga_i$ is a small linking
circle of $L_i$ then $[\ga_1],\dots,[\ga_r]$ is a basis for
$H_1(Y;\z)$. The dual basis for $H_2(Y;\z)$ can be represented by surfaces
$Z_1,\dots,Z_r$, where $Z_i$ is obtained by
capping off a Seifert surface of
$\ga_i\subset S^3$ by a disc. Here the Seifert surface should be
disjoint from the other components of $L$; this can be arranged by the obvious
tubing construction. It is clear that the bundle $P$ is
specified by the element of $(\z/2)^r$ whose $i$'th component
indicates whether $P|_{Z_i}$ is trivial or not. Without loss of
generality we may assume $P|_{Z_1}$ is non-trivial.

We form two other $\so$ bundles $P'\to Y'$
and $P''\to Y''$ as follows. The $3$-manifolds $Y'$ and $Y''$
are described by framed links $L'$ and $L''$ in $S^3$;
here $L'$ is obtained from $L$ by changing the framing of $L_r$
from $0$ to $-1$,
while $L''\subset L$ is the result of deleting
the component $L_r$. Then $b_1(Y')=b_1(Y'')=r-1$.
Let the bundles $P'$ and $P''$ both be specified by the element
of $(\z/2)^{r-1}$ which is the natural restriction of the vector
specifying $P$. By Floer's theorem we have an exact sequence
\[\hf^*(P')\oset{\al'}\to \hf^*(P)\oset\al\to \hf^*(P'').\]
Arguing as in the proof of Lemma~\ref{gen1_surg} we conclude that
the proposition holds for the bundle $P$, since by assumption it
holds for $P'$ and $P''$.

It remains to prove the proposition when $b_1(Y)=1$.
Again, let $Y$ be represented by a framed link $L$ as
in Lemma~\ref{repr}.
Choose a regular projection of $L$ and fix some
component $L_i$. It is well known that changing a
crossing within $L_i$ corresponds to $\pm1$ surgery on a
knot in $Y$ which bounds a surface of
genus $1$. By Lemma~\ref{gen1_surg} we may
therefore assume $L_i$ is unknotted. Since Kirby calculus allows
us to remove any unknotted component of $L$ with framing $\pm1$
(at the expense of
twisting the remainder of the link, but without changing framings or linking
numbers) we are left to consider the
manifold $S^1\times S^2$ described by the unknot in $S^3$ with
framing $0$. But in this case $\hf^*(P)=0$, so the proposition is
proved.\square

We now turn to homology spheres.

\be{thm}\label{fintyp2}
For any oriented homology $3$-sphere $Y$ there is a positive integer
$n$ such that
\[(u^2-64)^n=0\]
as an endomorphism of $\hfh^*(Y)$.
\end{thm}

\proof
Let $Y$ be an oriented homology
$3$-sphere and $\ga\subset Y$ a knot. Let $Y'$ and $Y''$ be the
result of $-1$ surgery and $0$ surgery on $\ga$, respectively.
There is then a long exact sequence
\be{equation}\label{ex1}
\dots\to \hf^{q+3}(Y'')\oset{\al''}\to \hf^q(Y)\oset\al\to \hf^q(Y')
\oset{\al'}\to \hf^q(Y'')\to\cdots.
\end{equation}
Note that $\del_0\al=\del_0$, since the cobordism from $Y$ to $Y'$
is negative definite and has no integral homology in dimension~$1$.
As explained in Section~\ref{umap} this sequence induces a
sequence of homomorphisms of $\q[u]$ modules
\be{equation}\label{ex2}
\hf^{q+3}(Y'')\to\hfh^q(Y)\to\hfh^q(Y')\to \hf^q(Y'').
\end{equation}
Moreover, it follows easily from Theorem~\ref{ucob2} and the exactness of
\Ref{ex1} that \Ref{ex2} is exact at the terms
$Y$ and $Y'$ for every $q$. We can then prove the theorem by the same
link reduction scheme as in the final paragraph of the proof of
Proposition~\ref{fintyp1a}, reducing the problem to $S^3$, where
it is trivial.\square

\be{thm}\label{uiso}
If $Y$ is any oriented homology $3$-sphere and $R$ any
associative ring in which $2$ is invertible then
$u:\hf^q(Y;R)\to \hf^{q+4}(Y;R)$ is an isomorphism for
$q\not\equiv4,5\mod8$.
\end{thm}

\proof By Theorem~\ref{fintyp2} the
$u$-map $\hf^q(Y;R)\to \hf^{q+4}(Y;R)$ is an
isomorphism for $q\equiv2,3\mod4$, since
in these degrees we have $\hf^q(Y;R)=\hfh^q(Y;R)$. Thus it only
remains to show that $u$ is an isomorphism for
$q\equiv0,1\mod8$. But this now follows from the exact sequence~\Ref{ex1}
and Lemma~\ref{gen1},
by the five-lemma and the same
induction scheme as we used in the proof of Theorem~\ref{fintyp2}.
\square

\be{cor}\label{mod4per}If $2$ is invertible in $R$ then
$\hfh^*(Y;R)$ and $\hf^*(Y;R)$ are both mod~$4$ periodic.
\end{cor}

\proof Recall that the Floer groups
are $\z/8$ graded. In the case of $\hfh^*(Y;R)$, it follows from
Theorem~\ref{fintyp2} that $u^2$ is invertible, hence
$u$ is an isomorphism of $\hfh^*(Y;R)$ onto itself.
As for $\hf^*(Y;R)$, note that for every $i$, the $u$-map is an isomorphism
$\hf^q(Y;R)\to\hf^{q+4}(Y;R)$ for either $q=i$ or $q=i-4$ (or both), by
Theorem~\ref{uiso}, hence $\hf^i(Y;R)$ and $\hf^{i+4}(Y;R)$ are
isomorphic.
\square

\section{Finite type of Donaldson invariants}\label{fintypsec}

We will now show that the nilpotency of $u^2-64$ on the reduced Floer
groups leads to a proof of the finite type conjecture of Kronheimer and
Mrowka in the simply-connected case.

\be{thm}\label{fintyp3}Let $X$ be a smooth, compact,
oriented $4$-manifold such that
$b_1(X)=0$ and $b_2^+(X)$ is odd.
Suppose there exists a splitting of $X$
along an embedded homology $3$-sphere $Y$,
\[X=X_1\cup_YX_2,\]
where $b_2^+(X_j)>0$ for $j=1,2$.
Then there exists a positive integer $n$ such
that for any homology orientation of $X$ and any $w\in H^2(X;\z)$
the Donaldson invariant
$D^w_X:\ba(X)\to\q$ satisfies
\[D^w_X((x^2-4)^nz)=0\]
for every $z\in\ba(X)$, where $x\in H_0$ is the point class.
In other words, $X$ has finite type. 
\end{thm}
Here  
\[\ba(X)=\text{Sym}(H_0(X;\q)\oplus H_2(X;\q)),\]
and $D^w_X$ is defined in terms of $\U2$ bundles $E$ over $X$ with $c_1(E)=w$.

\proof Because of the simplest blow-up formula \cite{Kotsch1} there is no
loss of generality in assuming $w_j=w|_{X_j}$ is not divisible by $2$ in
$H^2(X_j;\z)/\torsion$. To define relative invariants
\[D_j=D^{w_j}_{X_j}:\ba(X_j)\to\hf^*(\prtl X_j),\]
fix a metric on $Y$, and choose a generic
metric on $X_j\cup(\R_+\times\prtl X_j)$ which restricts to the given
product metric on the end. If $b_2^+(X_j)=1$ then $D_j$
depends on the chamber of the metric on $X_j$,
but this will not be reflected in our
notation. Donaldson's gluing theorem now says that for any
$z_j\in\ba(X_j)$ we have
\[D(z_1z_2)=D_1(z_1)\cdot D_2(z_2),\]
where $D=D^w_X$ and we use the natural pairing
$\hf^*(Y)\otimes\hf_*(Y)\to\q$ together with the identification
$\hf^*(\oline Y)=\hf_{5-*}(Y)$.

The crucial observation now is that $\del_nD_j=0$ for every $n\ge0$
because of the absence of reducible connections
(see section~\ref{red-conn}). Hence $D_j$ defines an invariant
\[\hat D_j:\ba(X_j)\to\hfh^*(\prtl X_j),\]
and we can use the natural pairing
$\hfh^*(Y)\otimes\hfh^{5-*}(\oline Y)\to\q$ in the gluing theorem:
\[D(z_1z_2)=\hat D_1(z_1)\cdot\hat D_2(z_2).\]
By Theorem~\ref{fintyp2} there exists a positive integer $n$ such
that $(u^2-64)^n=0$ on $\hfh^*(Y)$. This gives
\[\hat D_1((x^2-4)^nz_1)=((\frac u4)^2-4)^n\hat D_1(z)=0\]
for every $z_1\in\ba(X_1)$. (As in the proof of Proposition~\ref{red_count}
we can avoid twisted reducibles, since $b_1(X_1)=0$.)
But any class $z\in\ba(X)$ can be expressed as
$z=z_1z_2$ for some $z_j\in\ba(X_j)$, hence
\[D((x^2-4)^nz)=\hat D_1((x^2-4)^nz_1)\cdot\hat D_2(z_2)=0.\square\]

\be{thm}Let $X$ be a smooth, compact, simply-connected, oriented
$4$-manifold with $b_2^+(X)$ odd and $\ge3$. Then 
there exists a splitting of $X$ as in Theorem~\ref{fintyp3}, hence
$X$ has finite type.
\end{thm}

\proof By the classification of indefinite forms and Donaldson's theorem
we can express the intersection form of $X$ as an orthogonal sum
\[H_2(X;\z)/\torsion=V_1\oplus V_2\]
where both $V_1$ and $V_2$ contain vectors of positive square.
Since $X$ is simply-connected
we can invoke a theorem of Freedman
and Taylor \cite{FT}, \cite{Stong} which says that any orthogonal
splitting of the intersection form of $X$ is realized by some 
splitting of $X$ along an embedded homology sphere $Y$:
\[X=X_1\cup_YX_2.\square\]

\section{The $h$-invariant}\label{h-inv}

In this section we consider Floer groups with rational
coefficients, unless otherwise stated. We
define the $h$-invariant and establish two basic
properties: additivity under connected sums, and monotonicity with respect
to negative definite cobordisms.

\be{defn}For any oriented homology $3$-sphere define
\[h(Y)=\frac12(\chi(\hf^*(Y))-\chi(\hfh^*(Y))),\]
where $\chi$ is the Euler characteristic over $\q$.
\end{defn}

We will see in a moment that $h(Y)$ is always an integer.
Of course, by Taubes' theorem \cite{T2},
\[\chi(\hf^*(Y))=-2\lla(Y),\]
where $\lla$ is Casson's invariant.

Notice that we can identify $\del_n(Y)$
with $\del'_n(\oline Y)$ under the canonical isomorphisms
\[\hf^{5-q}(\oline Y)=\hf_q(Y)=(\hf^q(Y))^*.\]
Now let $B_*\subset \hf_*(Y)$ be the linear span of the
classes $\del_n$, $n\ge0$. We usually think of $\hf_q$ as the dual space
of $\hf^q$. Since
$\del_{2n+1}=\del_{2n}u$, and $u:\hf^0(Y)\to \hf^4(Y)$ is an isomorphism
by Theorem~\ref{uiso}, it follows that $\dim\,B_0=\dim\,B_4$ and
\[h(Y)=\dim\,B_4(Y)-\dim\,B_4(\oline Y).\]
As observed in Section~\ref{umap}, either $B_4(Y)=0$ or $B_4(\oline Y)=0$.

\be{prop}\label{h-char}
If $n$ is a non-negative integer then $h(Y)>n$ if and only if 
there exists an $x\in \hf^4(Y)$ such that $\del u^{2j}x=0$ for $0\le j<n$ 
but $\del u^{2n}x\neq0$.
\end{prop}
\proof 
It follows from Theorem~\ref{chain-map} that if $\del_{2k}$ lies in the
linear span of $\{\del_{2j}\}_{0\le j<k}$ then so does $\del_{2k+2}$.
Therefore, $h(Y)>n$ if and only if $\{\del_{2j}\}_{0\le j\le n}$ are
linearly independent. The proposition now follows because
\[\del_{2n}=\del u^{2n}\quad\text{on}\quad\bigcap_{j=0}^{n-1}\ker(\del_{2j}).
\square\]

\be{thm}\label{h_add}
$h(Y_1\#Y_2)=h(Y_1)+h(Y_2)$.
\end{thm}

We begin the proof of the theorem with five lemmas.

\be{lemma}\label{subadd}
If $h(Y_i)>0$ for $i=1,2$ then
$h(Y_1\#Y_2)\ge h(Y_1)+h(Y_2)$.
\end{lemma}
{\em Proof.} Let $W$ be the standard homology cobordism from
$\oline Y_1\cup\oline Y_2$ to $Y_1\#Y_2$, and let
\[W^*:\hf^p(Y_1)\otimes \hf^q(Y_2)\to \hf^{p+q}(Y_1\#Y_2)\]
be the homomorphism defined by $W$. We also consider the homomorphism
\[W^*_u:\hf^p(Y_1)\otimes \hf^q(Y_2)\to \hf^{p+q-4}(Y_1\#Y_2)\]
defined by cutting down moduli spaces over $W$ by four times the
$\mu$-class of a point. Then the proof of Theorem~\ref{chain-map} can be
adapted to show that
\[\del W^*_u(a_1\otimes a_2)=\pm2(\del a_1)(\del a_2).\]
On the other hand, moving the base point along a path
$[0,\infty)\times\{x_0\}$ in $\R_+\times(Y_1\#Y_2)$ one finds that
\[W^*_u(a_1\otimes a_2)=u W^*(a_1\otimes a_2).\]
Furthermore, one has
\[W^*(ua_1\otimes a_2)=uW^*(a_1\otimes a_2)\]
whenever $\del a_1=0$,
and similarly with the roles of $Y_1$ and $Y_2$ interchanged.
Now set $k_i=h(Y_i)$. By Proposition~\ref{h-char} and
Theorem~\ref{uiso} there is an element
$a_i\in\hf^0(Y_i)$ such that
\[\del u^{2r-1}a_i=\be{cases}
0 & \text{if $1\le r<k_i$}\\
\neq0 & \text{if $r=k_i$.}
\end{cases}\]
So if $1\le r_i\le k_i$ one has
\be{align*}
\del u^{2(r_1+r_2)-1}W^*(a_1\otimes a_2)
&=\del uW^*(u^{2r_1-1}a_1\otimes u^{2r_2-1}a_2)\\
&=\pm2(\del u^{2r_1-1}a_1)(\del u^{2r_2-1}a_2),
\end{align*}
and the lemma follows. \square

\be{lemma}Let $W$ be a smooth, compact, oriented $4$-manifold with boundary
$\prtl W=\oline{Y}_1\cup Y_2$, where both $Y_i$ are homology spheres.
If the intersection form of $W$ is negative definite and $H_1(W;\z)=0$ then
\[h(Y_2)\ge h(Y_1).\]
\end{lemma}

\proof This follows from Theorem~\ref{del-cob}.\square

\be{lemma}\label{posneg}
If $h(Y_1)=h(Y_2)>0$ then $h(Y_1\#\oline Y_2)=0$.
\end{lemma}
\proof If $h(Y_1\#\oline Y_2)>0$ then by the homology cobordism invariance
of $h$ one would have
\[h(Y_1)\ge h(Y_2)+h(Y_1\#\oline Y_2)>h(Y_2).\]
A similar argument applies if $h(Y_1\#\oline Y_2)<0$, since
$h(\oline Y)=-h(Y)$.\square

\be{lemma}\label{h_surg_weak}
Let $Y$ be an oriented homology $3$-sphere such that $h(Y)\ge0$,
and let $Y'$ be the result of
$-1$ surgery on a knot in $Y$ of genus $1$. Then
\[0\le h(Y')-h(Y)\le1.\]
\end{lemma}
Once we have established additivity of $h$ it will
be clear that the lemma holds without the assumption
$h(Y)\ge0$.

\proof Since the surgery cobordism $W$ from $Y$ to $Y'$ is negative definite
and satisfies $H_1(W;\z)=0$,
we have $h(Y)\le h(Y')$. Now set $n=h(Y)$.
To prove $h(Y')\le n+1$ we use the exact sequence~\Ref{ex1}.
Suppose $x\in\hf^4(Y')$ satisfies $\del u^{2j}x=0$ for $0\le j\le n$.
As in the proof of Lemma~\ref{gen1_surg} we have
\[\al'((u^2-64)x)=(u^2-64)\al'x=0,\]
so $(u^2-64)x=\al y$ for some $y\in\hf^4(Y)$. Since $\al u=u\al$
on $\ker\del_0$ by Theorem~\ref{ucob} we find that
for $0\le j\le n$,
\[\del u^{2j}y=\del\al u^{2j}y=\del u^{2j}(u^2-64)x=\del u^{2(j+1)}x.\]
Therefore, $0=\del u^{2n}y=\del u^{2(n+1)}x$. It follows that $h(Y')\le n+1$.
\square

Let $nY$ denote the $n$-fold connected sum
$\#_nY$ for $n\ge0$ (if $n=0$ we agree that $nY=S^3$),
and set $(-n)Y=n\oline Y$.
Let $S$ be the Brieskorn sphere $\Si(2,3,5)$.

\be{lemma}\label{h_poinc}
For any integer $n$ we have $h(n S)=n$.
\end{lemma}
\proof We may assume $n>0$. By Proposition~\ref{poinc_cup}~(i) we have
$h(S)=1$. Hence
$h(n S)\ge n$ by Lemma~\ref{subadd}. But $S$ is also the result
of $-1$ surgery on the negative $(2,3)$ torus knot, which has genus~1,
so $n$ applications of Lemma~\ref{h_surg_weak} gives $h(n S)\le n$.
Thus $h(n S)=n$.\square

{\em Proof of Theorem~\ref{h_add}.} Let $W$ be a smooth, compact, oriented,
connected $4$-manifold with boundary components $\oline Z_1,Z_2,\oline V_1,
\dots,\oline V_r$, where each component is a homology sphere. Suppose
$H_j(W;\z)=0$ for $j=1,2$ and $h(V_i)=0$ for each $i$.
We will show that $h(Z_1)=h(Z_2)$. If $h(Z_1)$ and $h(Z_2)$ are both zero
then there is nothing to prove, so after perhaps reversing orientations
we may assume $h(Z_2)\le h(Z_1)>0$.

Since $h(V_i)=0$ we can find $\rho_i\in CF^0(V_i)$ such that
$d\rho_i=\del'$. Let $\hat W=W\cup(\R_+\times\prtl W)$ have a tubular
end metric. Then $0$-dimensional moduli spaces over $\hat W$ with
the chain $\theta+\rho_i$ as ``flat limit'' over the end $\R_-\times V_i$
define a degree preserving homomorphism
\[f:\hf^*(Z_1)\to\hf^*(Z_2)\]
which satifies $fu=uf$ on $\ker \del_0$ and $\del_0f=\del_0$.
This implies $h(Z_2)\ge h(Z_1)$, so $h(Z_1)=h(Z_2)$.

To prove the theorem, set $k_i=h(Y_i)$,
$k=k_1+k_2$, and let $W$ have boundary components
$Y_1\#Y_2$, $k\oline S$, $\oline Y_1\#k_1 S$, and
$\oline Y_2\#k_2 S$.\square

\be{thm}\label{h_cob}
Let $W$ be a smooth, compact, oriented $4$-manifold with boundary
$\prtl W=\oline{Y}_1\cup Y_2$, where both $Y_i$ are homology spheres.
Suppose the intersection form of $W$ is negative definite.
Then
\[h(Y_2)\ge h(Y_1),\]
with strict inequality if the intersection form is not diagonal over
the integers.
\end{thm}

\proof Let $L$ be the intersection form of $W$.
Then $\ti Y=Y_2\#\oline Y_1\#S$ bounds a smooth, compact,
oriented $4$-manifold with negative definite intersection form $L\oplus(-E_8)$,
which is not diagonal over the integers. Hence $h(\ti Y)\ge1$ by
Corollary~\ref{ex-int}. Since $h(S)=1$ and $h$ is additive, we deduce
$h(Y_2)\ge h(Y_1)$.

(It is possible that a more direct proof of this can be found by first
surgering away the free part of $H_1(W;\z)$ and then analyzing
the abelian flat $\su$ connections over $W$ as in \cite{D2},
Section~4b.)

If $L$ is not diagonal then $h(Y_2)-h(Y_1)=h(Y_2\#\oline Y_1)>0$
because $Y_2\#\oline Y_1$ bounds a smooth, compact,
oriented $4$-manifold with the same intersection form as $W$. 
\square

\nocite{Fr1}

\textsc{Harvard University, Cambridge, MA 02138}
\end{document}